\documentclass[leqno, 12pt]{amsart}
\usepackage{amscd}
\usepackage{amssymb}
\usepackage[all, poly]{xy}

\numberwithin{equation}{subsection}

\newcommand{\nc}{\newcommand}
\newcommand{\rc}{\renewcommand}

\nc{\lab}{      \label}
\nc{\npp}{{     \newpage\setcounter{page}{0}    }}
\nc{\setpart}{{         \setcounter{part}       }}
\nc{\setpage}{{         \setcounter{page}       }}
\nc{\setsection}{{      \setcounter{section}    }}
\nc{\nd}{ $$\text{ This version is preliminary and approximate, 
it is not for distribution. }$$ }

\nc{\noi}{{\noindent}}
\nc{\nop}{{\noindent {\bf Proof.}} }

\nc{\indd}{{ ${} \ \ \ \ \  \ \        {} $}}   

\nc{\bss}{{\backslash}}           

\nc{\barr}{     \overline       }               
\nc{\ud}{       \underline      }               

\nc{\sub}{{     \subseteq       }}         
\nc{\subb}{{    \supseteq       }}         
\nc{\nsub}{{    \nsubseteq      }}         
\nc{\nsubb}{{   \nsupseteq      }}         %

\nc{\nin}{{     \notin  }}      

\nc{\ti}{\tilde}              
\nc{\wtl}{\widetilde}         

\nc{\hatt}{\widehat}                            
\nc{\hata}{{    \bbb{ \hat{} }          }}      

\nc{\ch}{\check}                                
\nc{\cha}{{     \bbb{ \check{} }        }}      

\nc{\lb}{\langle}                                       
\nc{\rb}{\rangle}

\nc{\lB}{       \left(  }                               
\nc{\rB}{       \right) }
\nc{\BBl}{{     \bbb{ \left( \right.}   }}              
\nc{\BBr}{{     \bbb{ \left. \right)}   }}

\nc{\h}{{       \hslash }}      

\nc{\All}{{     \forall }}
\nc{\Ex}{{      \exists         }}

\nc{\yy}{\infty}

\nc{\se}{       \section                }
\nc{\sus}{      \subsection             }
\nc{\sss}{      \subsubsection          }

\nc{\Lemm}{     \subsection{Lemma}              }
\nc{\lemm}{     \subsubsection{Lemma}           }
\nc{\slemm}{    \subsubsection*{Lemma}          }

\nc{\Pro}{      \subsection{Proposition}        }
\nc{\pro}{      \subsubsection{Proposition}     }
\nc{\spro}{     \subsubsection*{Proposition}    }

\nc{\Corr}{     \subsection{Corollary}          }
\nc{\corr}{     \subsubsection{Corollary}       }
\nc{\scorr}{    \subsubsection*{Corollary}      }

\nc{\Theo}{     \subsection{Theorem}            }
\nc{\theo}{     \subsubsection{\bf Theorem}         }
\nc{\stheo}{    \subsubsection*{\bf Theorem}                }

\nc{\remm}{     \subsubsection{Remarks}         }
\nc{\sremm}{    \subsubsection*{Remarks}                }
\nc{\rema}{     \subsubsection{Remarks}         }

\nc{\Conj}{     \subsection{\bf Conjecture} }
\nc{\conj}{     \subsubsection{\bf Conjecture}              }
\nc{\sconj}{    \subsubsection*{\bf Conjecture}             }

\nc{\ex}{       \subsubsection{Example}         }
\nc{\sex}{      \subsubsection*{Example}                }
\nc{\exs}{      \subsubsection{Examples}                }
\nc{\sexs}{     \subsubsection*{Examples}               }

\rc{\AA}{{\mathcal A}} 
\nc{\CC}{{\mathcal C}}
\nc{\DD}{{\mathcal D}}
\nc{\EE}{{\mathcal E}}
\nc{\FF}{{\mathcal F}}
\nc{\GG}{{\mathcal G}}
\nc{\HH}{{\mathcal H}}
\nc{\II}{{\mathcal I}}
\nc{\JJ}{{\mathcal J}}
\nc{\KK}{{\mathcal K}}
\nc{\LL}{{\mathcal L}}
\nc{\MM}{{      \mathcal M        }}
\nc{\NN}{{\mathcal N}}
\nc{\OO}{{\mathcal O}}
\nc{\PP}{{\mathcal P}}
\nc{\QQ}{{\mathcal Q}}
\nc{\RR}{{\mathcal R}}
\rc{\SS}{{\mathcal S}}
\nc{\TT}{{\mathcal T}}
\nc{\UU}{{\mathcal U}}
\nc{\VV}{{\mathcal V}}
\nc{\WW}{{\mathcal W}}
\nc{\ZZ}{{\mathcal Z}}
\nc{\XX}{{\mathcal X}}
\nc{\YY}{{\mathcal Y}}

\nc{\A}{{\mathbb A }}
\nc{\B}{{\mathbb B}}
\nc{\C}{{\mathbb C}}
                \nc{\cc}{{\mathbb C}}
\nc{\Cs}{{\mathbb C^*}}
                \nc{\cs}{{\mathbb C^*}}
                \nc{\ccs}{{\mathbb C^*}}
\nc{\D}{{\mathbb D}}
\nc{\E}{{\mathbb E}}
\nc{\F}{{\mathbb F}}
\nc{\G}{{\mathbb G}}
        \nc{\hH}{{\mathbb H}}
\nc{\I}{{\mathbb I}}
\nc{\J}{{\mathbb J}}
\nc{\M}{{\mathbb M}}
\nc{\N}{{\mathbb N}}
        \nc{\oO}{{\mathbb O}}
        \nc{\pP}{{\mathbb P}}      
\nc{\Q}{{\mathbb Q}}
\nc{\R}{{\mathbb R}}
        \nc{\sS}{{\mathbb S}}
\nc{\T}{{\mathbb T}}
\nc{\U}{{\mathbb U}}
\nc{\V}{{\mathbb V}}
\nc{\W}{{\mathbb W}}
\nc{\Z}{{\mathbb Z}}
\nc{\X}{{\mathbb X}}
\nc{\Y}{{\mathbb Y}}

\let\H\hH

\let\P\pP

\nc{\fA}{{\mathfrak A}}
\nc{\fB}{{\mathfrak B}}
\nc{\fC}{{\mathfrak C}}
\nc{\fD}{{\mathfrak D}}
\nc{\fE}{{\mathfrak E}}
\nc{\fF}{{\mathfrak F}}
\nc{\fG}{{\mathfrak G}}
\nc{\fH}{{\mathfrak H}}
\nc{\fI}{{\mathfrak I}}
\nc{\fJ}{{\mathfrak J}}
\nc{\fK}{{\mathfrak K}}
\nc{\fL}{{\mathfrak L}}
\nc{\fM}{{\mathfrak M}}
\nc{\fN}{{\mathfrak N}}
\nc{\fO}{{\mathfrak O}}
\nc{\fP}{{\mathfrak P}}
\nc{\fQ}{{\mathfrak Q}}
\nc{\fR}{{\mathfrak R}}
\nc{\fS}{{\mathfrak S}}
\nc{\fT}{{\mathfrak T}}
\nc{\fU}{{\mathfrak U}}
\nc{\fV}{{\mathfrak V}}
\nc{\fW}{{\mathfrak W}}
\nc{\fZ}{{\mathfrak Z}}
\nc{\fX}{{\mathfrak X}}
\nc{\fY}{{\mathfrak Y}}
\nc{\fa}{{\mathfrak a}}
\nc{\fb}{{\mathfrak b}}
\nc{\fc}{{\mathfrak c}}
\nc{\fd}{{\mathfrak d}}
\nc{\fe}{{\mathfrak e}}
\nc{\ff}{{\mathfrak f}}
\nc{\fg}{{\mathfrak g}}
\nc{\fh}{{\mathfrak h}}
\nc{\fiI}{{\mathfrak i}}  
        \nc{\ffi}{{\mathfrak i}}  
\nc{\fj}{{\mathfrak j}}
\nc{\fk}{{\mathfrak k}}
\nc{\fl}{{\mathfrak{l}}}
\nc{\fm}{{\mathfrak m}}
\nc{\fn}{{\mathfrak n}}
\nc{\fo}{{\mathfrak o}}
\nc{\fp}{{\mathfrak p}}
\nc{\fq}{{\mathfrak q}}
\nc{\fr}{{\mathfrak r}}
\nc{\fs}{{\mathfrak s}}
\nc{\ft}{{\mathfrak t}}
\nc{\fu}{{\mathfrak u}}
\nc{\fv}{{\mathfrak v}}
\nc{\fw}{{\mathfrak w}}
\nc{\fz}{{\mathfrak z}}
\nc{\fx}{{\mathfrak x}}
\nc{\fy}{{\mathfrak y}}

\nc{\al}{{\alpha }}
\nc{\be}{{\beta }}
\nc{\ga}{{\gamma }}
\nc{\de}{{\delta }}
\nc{\del}{{\partial }}
\nc{\ep}{{\varepsilon }}
\nc{\vap}{{\epsilon }}

\nc{\ze}{{\zeta }}
\nc{\et}{{\eta }}
\rc{\th}{{\theta }}
\nc{\vth}{{\vartheta }}

\nc{\io}{{\iota }}
\nc{\ka}{{\kappa }}
\nc{\la}{{\lambda }}
\nc{\vrho}{{\varrho}}
\nc{\si}{{\sigma }}
\nc{\ups}{{\upsilon }}
\nc{\vphi}{{\varphi }}
\nc{\om}{{\omega }}

\nc{\Ga}{{\Gamma }}
\nc{\De}{{\Delta }}
\nc{\nab}{{\nabla}}
\nc{\Th}{{\Theta }}
\nc{\La}{{\Lambda }}
\nc{\Si}{{\Sigma }}
\nc{\Ups}{{\Upsilon }}
\nc{\Om}{{\Omega }}

\nc{\toc}{{\tableofcontents}}
\nc{\addl}{     \addcontentsline{toc}{subsection}       }

\nc{\GK}{{      G(\KK)          }}
\nc{\GO}{{      G(\OO)          }}

\nc{\BK}{{  B(\mathcal K)               }}
\nc{\BO}{{  B(\mathcal O)               }}
\nc{\BKz}{{  B(\mathcal K)_0    }}
\nc{\NK}{{  N(\mathcal K)               }}
\nc{\NO}{{  N(\mathcal O)               }}
\nc{\TK}{{  B(\mathcal K)               }}
\nc{\TO}{{  T(\mathcal O)               }}
\nc{\LO}{{  L(\mathcal O)               }}

\nc{\htt}{  \text{ht}}

\nc{\pt}{  \text{pt}}

\newtheorem*{theorem}{Theorem}

\newtheorem*{ta}{Theorem A}
\newtheorem*{TheTheorem}{The $d^2=0$ Theorem}

\nc{\BFA}{{\bf A}}
\nc{\BFB}{{\bf B}}
\nc{\BFE}{{\bf E}}
\nc{\BFF}{{\bf F}}
\nc{\BFG}{{\bf G}}
\nc{\BFV}{{\bf V}}

\nc{\Ext}{\operatorname{Ext}}
\nc{\Hom}{\operatorname{Hom}}
\nc{\End}{\operatorname{End}}
\nc{\Ker}{\operatorname{Ker}}
\nc{\Imm}{\operatorname{Im}}
\nc{\triv}{\operatorname{triv}}

\nc{\bfb}{{\bf b}}

\nc{\Id}{\operatorname{Id}}
\nc{\opp}{\operatorname{opp}}

\nc{\ic}{\operatorname{IC}}
\nc{\mof}{\operatorname{mof-}}
\nc{\pic}{{\mathbb P}\operatorname{IC}}
\nc{\bic}{{\bf ic}}

\nc{\ICS}{{\bf IC}}

\nc{\BFP}{{\bf P}}

\nc{\Vect}{\operatorname{Vec}}
\nc{\mixed}{\operatorname{mixed}}
\nc{\OPH}{\operatorname{H}}

\nc{\VP}{{\mathbb V}P}

\nc{\Rel}{\operatorname{Rel}}
\nc{\Lie}{\operatorname{Lie}}
\nc{\spann}{\operatorname{span}}
\nc{\ev}{\operatorname{ev}}

\long\def\MSC#1\EndMSC{\def\arg{#1}\ifx\arg\empty\relax\else
     {\par\narrower\noindent%
     2000 Mathematics Subject Classification: #1\par}\fi}

\begin{document}

\title[]{
Perverse sheaves, Koszul $\ic$-modules, \\
and the quiver for the category $\OO$
}

\author{       Maxim Vybornov                }
\address{Dept. of Mathematics, MIT, 
77 Mass Ave, Cambridge MA 02139, USA}
\email{                vybornov@math.mit.edu }

\curraddr{22 Rockland St, Newton, MA 02458, USA}

\begin{abstract} 
For a stratified topological space
we introduce the category of $\ic$-modules, which are
linear algebra devices
with the relations described by the equation $d^2=0$.
We prove that the category of (mixed) $\ic$-modules
is equivalent to the category of (mixed)
perverse sheaves for flag varieties. As an application, 
we describe an algorithm 
calculating the quiver underlying the BGG category 
$\mathcal O$
for arbitrary simple Lie algebra, thus answering
a question which goes back to I.~ M.~ Gelfand.
\end{abstract}



\maketitle


\begin{center}
\emph{Dedicated to George Lusztig on the occasion of his 60-th birthday}
\end{center}

\se{Introduction}

\sus{} Motivated by R.~ MacPherson's cellular perverse sheaves,
W.~ Soergel's studies of the Bernstein-Gelfand-Gelfand (BGG)
category $\mathcal O$,
and 
L.~ Saper's 
$L$-modules we introduce here, in Section \ref{DefSec}, 
an abelian category
of \emph{$\ic$-modules} designed to be a 
derived-category-free
model
of the category of constructible perverse sheaves 
if the latter is Koszul (cf. \ref{AssumpKoszul}
for the remarks beyond the Koszul case). 
In Section \ref{SecFlagProof} we prove the following.

\begin{ta} The category of Schubert-constructible 
perverse sheaves on a flag variety is equivalent to the category
of $\ic$-modules.
\end{ta}

Moreover, we suggest a relationship between perverse 
sheaves and $\ic$-modules 
for a much wider class of stratified spaces.

\sus{} There exists a number of derived-category-free models of 
perverse sheaves. 
These include $\DD$-modules, cf. e.g. \cite{KS}
(which actually predate perverse sheaves),
gluing constructions \cite{Be86, MV86}, 
and ``quiver'' presentations, cf. e.g. \cite{MV88, GMV, BrGr, Br}.
We believe that our ``elementary'' construction has the advantage
of being well suited for \emph{mixed} perverse sheaves.
We construct
a natural functor {f}rom mixed perverse sheaves to 
mixed $\ic$-modules which is an equivalence
in our two main examples: flag 
varieties and simplicial complexes, 
cf. Section \ref{WeightSec}.

\sus{} Let $G$ be a complex simple algebraic group,
$\fg=\Lie G$ its Lie algebra, $B$ a Borel
subgroup, and $\fb=\Lie B$ its Lie algebra.
Let $\fh$ be a Cartan of $\fb$, $\Phi$ the root system, and let 
$W$ be the Weyl group. Let $\OO_0$ be a regular block of the BGG
category $\OO$, cf. \cite{BGG71, BGG73, BGG76},
equivalent to Schubert-stratification constructible perverse sheaves
on $G/B$, \cite{BB, BBD}. 
The simple objects of $\OO_0$ are the simple 
$\fg$-modules 
with the highest weights $w(\rho)-\rho$ where 
$\rho=\frac{1}{2}\sum_{\al\in \Phi^+} \al$.
If $\fg=\mathfrak{sl}_2$, then $\OO_0$
is equivalent to the following category of linear 
algebra data
\begin{equation}\nonumber
\qquad   \xymatrix{
        V_0 \ar@/_/[r]_{b} &  
        V_1  \ar@/_/[l]_{a} & 
}
\end{equation} 
with the relation $b\circ a=0$. This is the 
{\em quiver} of the category $\OO_0$
for $\mathfrak{sl}_2$. The question of finding 
such a quiver for arbitrary simple 
Lie algebra $\fg$
goes back to
I.~ M.~ Gelfand, cf. \cite{Gelf}. 
Seminal breakthroughs with implications
for this question had been made in particular in
\cite{KL79, BB, So90, BGS}.

It is well known that the quiver has
vertices labeled by the elements of the Weyl group $W$.
If $y, w\in W$, then 
for every arrow {f}rom $y$ to $w$ there is an arrow {f}rom 
$w$ to $y$. Moreover, the number of arrows from $y$ to $w$
is equal to $\mu(y, w)$ -- the coefficient of the power
$\frac{1}{2} (l(w)-l(y)-1)$ of the Kazhdan-Lusztig polynomial 
\cite{KL79}. However, the relations between arrows 
were not known explicitly
and were difficult to obtain algorithmically, except
for a few small rank cases, cf. \cite{I82, I85, St03}.

In Section \ref{Algm} we use $\ic$-modules 
to encode the relations in the form $d^2=0$ and describe
a \emph{computer-executable algorithm} to calculate them 
explicitly. The algorithm is a refinement of an algorithm of 
C. ~ Stroppel \cite{St03}. It
works for arbitrary simple 
Lie algebra $\fg$, thus settling the Gelfand's question 
for the category $\mathcal O$. A detailed implementation of
the algorithm for $\fg = \mathfrak{sl}_3$ 
is worked out in the Appendix (Section 8).

\sus{} We believe that $\ic$-modules should also be useful for
the study of Harish-Chandra modules and Soergel's 
conjecture, cf. \cite[Basic conjecture 1]{So01}, 
as well as
toric varieties, affine Grassmannians, and 
locally symmetric spaces,
cf. Section \ref{FinalRemarks}.

\sus*{Acknowledgment} 
I would like to thank Robert MacPherson, 
Ivan Mirkovi\' c and Wolfgang Soergel
for some crucial suggestions and observations.
I am grateful to the anonymous referee as well as 
A.~ Beilinson, J.~ Bernstein, T.~ Braden,
G.~ Faltings, O.~ Gabber, S.~ Gelfand, G.~ Lusztig, V.~ Ostrik,
A.~ Polishchuk,  A.~ Postnikov, C. ~ Stroppel, and
K.~ Vilonen
for very useful discussions,
and to MPIfM, IH\' ES, 
and Universit\" at Freiburg
for their 
hospitality and support. 
I am grateful to Y.~ Drozd, S.~ Ovsienko, and V.~ Futorny
who introduced me to Lie algebras and quivers during my college years. 
The research was supported in part by the NSF Postdoctoral 
Research Fellowship.

\se{Definitions and The $d^2=0$ Theorem}\label{DefSec}

\sus{Setup}
Our sheaves are assumed to be sheaves of finitely generated 
$R$-modules over a ring $R$. 
Unless indicated otherwise, $R$ is assumed to be 
the field $\C$, 
but many results are
true for other ground rings. 

Let $X$ be a topological space equipped with a sheaf of rings
and let
$$
X=\bigsqcup_{S\in\SS} S
$$
be a partition of $X$ into a finite disjoint union
of locally closed strata. Let $\dim: \SS \to \Z_{\geq 0}$
be a dimension function such that $\dim S\leq \dim S'$
if $\barr S\subseteq \barr S'$ and let us fix 
a perversity function
$p: \Z_{\geq 0}\to \Z$,
cf. \cite{GM, BBD}.

Let $L$ be a local system on $S$.
Consider 
the $\ic$-sheaf $\ICS(\barr S, L)=\ICS_p(\barr S, L)$, 
cf. \cite{GM80, GM, BBD}
corresponding to $L$ and the perversity $p$
on the closure $\barr S$
of the stratum $S$, and let 
$$
IH^*(\barr S, L)=IH^*_p(\barr S, L)= \H^*(\ICS_p(\barr S,L))
$$
be the hypercohomology of the sheaf $\ICS_p(\barr S, L)$.
{F}rom now on we will omit 
perversity from the notation. Perversity is always assumed to 
be the middle one for algebraic varieties with
algebraic stratifications.
Notice that $IH^*(\barr S, L)$ is naturally a graded
$H^{*}(X)$-module, where $H^{*}(X)$ is the cohomology of $X$
with coefficients in $R$. 

{F}rom now on we assume that all strata are \emph{contractible}
unless explicitly specified otherwise, e.g. in section 
\ref{LocalSystem}. 

\sus{Definition}\label{mixedDefinition}
 A \emph{mixed} $\ic$-module $V$ 
is the following data:
\begin{enumerate}
\item for every stratum $S\in \SS$,
a $\Z$-graded $R$-vector space (stalk) $\oplus_{i\in \Z}V^i_S$ on $S$,
\item for every pair of incident 
(i.e. one is in the closure of the other)
strata $S$ and $S'$ such that $S\neq S'$, 
a boundary map
$$
v(S,S')\in \Hom^{(1,1)}_{H^{*}(X)}(IH^{*}(\barr S, V^{\bullet}_S),
IH^{*}(\barr S', V^{\bullet}_{S'})),
$$
\end{enumerate}
where $\Hom^{(1,1)}$ is the space of degree $(1,1)$ morphisms of 
$H^{*}(X)$-modules
$IH^{*}(\barr S, V^{\bullet}_S)$ bigraded by the degree 
of $IH^{*}$ and the degree of $V^{\bullet}_S$.
The data is subject to the  
\begin{equation}\label{ChainComplexAxiom}
\text{\rm Chain Complex Axiom:}\qquad  d^2=0,
\end{equation}
where $d$ is the degree $(1,1)$ map
\begin{equation}
d=\bigoplus_{S,S'} v(S,S'): \oplus_S IH^{*}(\barr S, V^{\bullet}_S)\to
\oplus_S IH^{*}(\barr S, V^{\bullet}_S).
\end{equation}
The bigraded space $\oplus_S IH^{*}(\barr S, V^{\bullet}_S)$ with the 
differential 
$d$ is called the \emph{total complex}
of the $\ic$-module $V$.
Its cohomology is naturally bigraded.

\sus{}\label{nonMixedDefinition}
A \emph{non-mixed} $\ic$ module is defined 
exactly as above, except the 
stalks $V_S$ are just (non-graded) $R$-vector spaces, 
and boundary maps 
$$
v(S,S')\in \Hom^{1}_{H^{*}(X)}(IH^{*}(\barr S),
IH^{*}(\barr S'))\otimes \Hom_{R}(V_S,V_{S'}) 
$$
and the differential $d$ are maps of degree $1$. The cohomology 
of the total complex has only one natural grading. 

\sus{} A morphism of $\ic$ modules is the collection
of maps between stalks commuting with the boundary maps.
The abelian category of non-mixed $\ic$-modules of finite length
(i.e. having a finite filtration
with simple subquotients) will be denoted by $\AA_{\SS}(X)$. 

The category of mixed $\ic$-modules of finite length 
will be denoted by
$\AA^{\mixed}_{\SS}(X)$. It is a mixed category in the sense 
of \cite[4.1]{BGS}. 
A mixed $\ic$-module $V$ is called \emph{pure} of weight $m$ 
if it is concentrated in 
degree $-m$ i.e.,  $V_S^i=0$ unless $i=-m$ for all $S\in \SS$.
We have the obvious functor 
$v: \AA^{\mixed}_{\SS}(X)\to \AA_{\SS}(X)$ forgetting the mixed
structure.

\sus{Local Systems}\label{LocalSystem} One could further
generalize the definition of an $\ic$-module
to the case of non-simply-connected strata. In general,
the stalk $V_S$ is a local system on $S$, and the boundary 
maps are
$$
v(S,S')\in 
\Hom^1_{H^{*}(X)}(IH^{*}(\barr S, V_S),IH^{*}(\barr S', V_{S'})),
$$
where $\Hom^1$ are the degree $1$ morphisms of $H^{*}(X)$-modules.
(We consider here the non-mixed case to simplify the notation.)
The {\rm Chain Complex Axiom} and the total complex are precisely 
as above.

{F}rom now on we will consider non-mixed $\ic$-modules
unless explicitly specified otherwise, e.g. in 
Section \ref{WeightSec}.

\sus{Verdier duality} 
We will only consider here the middle perversity case in order 
to simplify the notation. For an $\ic$-module $V$ the stalks of 
its Verdier dual $\ic$-module $DV$
are defined as follows:
$$
DV_S=V^*_S,
$$
where $V^*_S$ is the dual vector space.
Notice that 
$$
(V_{S}\otimes IH^{i}(\barr S))^*=V^*_{S}\otimes IH^{-i}(\barr S)
$$
since according to \cite[5.3]{GM} we have
$
IH^i(\barr S)=(IH^{-i}(\barr S))^*
$.
Now the boundary maps
\begin{equation}\nonumber
\begin{matrix}
dv^i(S,S'): & V_{S}^*\otimes IH^{i}(\barr S)
& \to &  
V_{S'}^*\otimes IH^{i+1}(\barr S') \\
& || && ||  \\
& V_{S}^*\otimes (IH^{-i}(\barr S))^* & \to & V_{S'}^*\otimes 
(IH^{-i-1}(\barr S'))^*
\end{matrix}
\end{equation}
of $DV$ are $dv^i(S,S')=(v^{-i-1}(S', S))^*$, where 
$$ 
v^{-i-1}(S', S): V_{S'}\otimes IH^{-i-1}(\barr S')
\to V_{S}\otimes IH^{-i}(\barr S)
$$
are the boundary maps of $V$.
The Chain Complex Axiom (\ref{ChainComplexAxiom})
for $DV$ is obviously satisfied.
We have constructed a contravariant functor 
$D: \AA_{\SS}(X)\to \AA_{\SS}(X)$. Clearly, $D\circ D=\Id$.

\sus{}
Denote the category of $\SS$-constructible 
perverse sheaves of finite length on $X$ by $\PP_{\SS}(X)$, 
cf. \cite{BBD}. We present here a more precise formulation of 
Theorem A.
The \emph{raison d'\^ etre} of these notes is the following.

\begin{TheTheorem}\label{TheTheorem}
If $X$ is a flag variety stratified by Schubert cells,
or a simplicial complex stratified by simplices, then
\begin{itemize}
\item[{(i)}] the categories
$$
\PP_{\SS}(X)\simeq \AA_{\SS}(X)
$$
are equivalent,
\item[{(ii)}] the total complex
of an $\ic$-module calculates the hypercohomology
of the corresponding perverse sheaf. 
\end{itemize}
\end{TheTheorem}

\sss*{\bf Remark: The $d^2=0$ Thereom for 
simplicial complexes}\label{EvidenceSimp} 
For perverse sheaves on simplicial complexes
constructible with respect to the triangulation
we have
$$
IH^*(\barr\Delta)=\C[\delta(\dim\Delta)],
$$
where $\Delta$ is a simplex and $\delta$
is a cellular perversity, \cite{M93, P, V00}.
The $\ic$-modules in this case are 
\emph{simplicial perverse sheaves}
of R.~ MacPherson and
the statement of the The $d^2=0$ Thereom
is a theorem due to MacPherson \cite{M93}, 
cf. \cite{P, V98, V00} for an alternative proof.

\sus{Restrictions, relaxing remarks, and conjectures}
\label{assump} 

We use only a limited number of properties 
(some are listed below)
of the category
of perverse sheaves in order to prove the Theorem.
It is natural to conjecture that under these restrictions
The $d^2=0$ Theorem holds. (Note that
only the statement (ii) is in question.)

\sss{}\label{FiniteStratification} Assume that we have only finitely 
many strata $S\in\SS$. This assumption can be relaxed, but we will 
not pursue it here.

\sss{}\label{AssumpContractible} For the purposes of this text 
we will assume that the strata are contractible.
This assumption could be relaxed: The $d^2=0$ Theorem
with a local system version of $\ic$-modules
is a tautology for a stratification
with one stratum, for example.

\sss{}\label{EnoughProjectives} Assume that the category
$\PP_{\SS}(X)$ has enough projectives i.e., every object 
has a projective resolution.

\sss{} \label{AssumpKoszul} 
We assume that the category $\PP_{\SS}(X)$ is Koszul
i.e., it is equivalent to a category of finitely generated modules
over a Koszul algebra.
When the category $\PP_{\SS}(X)$ is not Koszul
one can still reconstruct it from the 
algebra $\Ext^{*}_{\PP_{\SS}(X)}(L,L)$,
where $L$ is the direct product of all simple objects in $\PP_{\SS}(X)$
with the canonical
$A_{\infty}$-structure, cf. e.g. \cite[Problem 2]{K}. 
One has to take
this $A_{\infty}$-structure into account 
in order to generalize the $\ic$-modules
beyond the Koszul case.
\sss{} \label{AssumpFaith}
For any two incident strata $S$ and $S'$ 
and two simple objects $\LL_1=\ICS(S, L_S)$ and 
$\LL_2=\ICS(S', L_{S'})$ 
in $\PP_{\SS}(X)$ we assume that we have the isomorphisms
$$
\Ext^i_{\PP_{\SS}(X)}(\LL_1, \LL_2) =
\Ext^i_{D(X)}(\LL_1, \LL_2) 
=
\Hom^{i}_{H^*(X)}(\H^*(\LL_1),\H^*(\LL_2)),
$$
and
$
\Ext^i_{\PP_{\SS}(X)}(\LL_1, \LL_2) = \Ext^i_{D(X)}(\LL_1, \LL_2) =0
$
for any two non-incident strata $S$ and $S'$.
Here $\Hom^i$
is the space of $H^{*}(X)$-morphisms
of degree $i$. 
In particular, the representation of the algebra 
$\Ext^{*}_{\PP_{\SS}(X)}(L,L)$ in the space
$\oplus_S IH^*(\barr S, L_S)$
is faithful.
One could generalize $\ic$-modules 
using another faithful representation of the 
algebra $\Ext^{*}_{\PP_{\SS}(X)}(L, L)$.

\se{Proof of Theorem A}\label{SecFlagProof}

In this section we prove The $d^2=0$ Theorem, 
and thus Theorem A,
following a suggestion of W.~ Soergel. 
In fact our proof
works for any stratified space with the restrictions
\ref{assump} and one additional condition: $R_X[\dim X]$ 
is a simple object of $\PP_{\SS}(X)$.

\sus{} Let $G$ be a simple algebraic group over $\C$,
$B\subseteq P\subseteq G$ be a Borel subgroup 
and a parabolic subgroup, $W_P$ the subgroup
of the Weyl group $W$
and let $X = G/P$ be the corresponding partial flag variety,
stratified by Schubert cells $X_w$, where $w$
is a coset in $W/W_P$ which we identify with the set
of strata $W_{\SS}$ of $X$.
Denote by 
$\PP_{\rm{Schubert}} (G/P)$ 
the category of Schubert-stratification
constructible perverse sheaves on $G/P$.

\sus{Theorem}\label{SoergelSatz} \cite{So90, G}
\emph{For any $y, w\in W_{\SS}$ the natural map
$$
\Ext^*_{D(X)}(\ICS_y, \ICS_w)=\Hom^{*}_{H^*(X)}
(IH^*(\barr X_y),IH^*(\barr X_w))
$$
is an isomorphism of graded spaces. (Here $\Hom^{i}$
denotes degree $i$ morphisms of $H^*(X)$-modules.)
Moreover, the left and right hand sides are 
isomorphic as graded algebras.}

\sus{Koszul complex}\label{KoszulComplexSection}
In this subsection $A=\oplus_{i\geq 0}A_i$ is an arbitrary Koszul 
algebra with semisimple $A_0=k$ over the field $\C$. 
All tensor products till the end of Section \ref{SecFlagProof}
are over $k$
unless specified otherwise.
Let $A^{!}$ be the quadratic dual
algebra and let
us consider the Koszul complex
\cite{BGS}
\begin{equation}\label{KoszulComplex}
\cdots A\otimes{}^{\ast}(A^{!}_2)\to  
A\otimes{}^{\ast}(A^{!}_1)\to A\to k.
\end{equation}
Recall the construction of the differential from \cite[2.7-2.8]{BGS}.
Identify 
$A\otimes {}^{\ast}(A^{!}_i)=\Hom_{-k}(A^{!}_i, A)$.
Let $\Id_{A_1}=\sum a'_{\al}\otimes a_{\al}\in A^{!}_1\otimes A_1$
under the canonical isomorphism 
$\Hom_{k}(A_1, A_1)=A^{!}_1\otimes A_1$.
Observe that $a_{\al}$ (resp. $a'_{\al}$) is a basis in $A_1$  
(resp. $A^{!}_1$). Then the differential
$$ 
d: A\otimes {}^{\ast}(A^{!}_{i+1})\to A\otimes {}^{\ast}(A^{!}_i)
$$
is constructed as follows: 
$(df)(a)=\sum f(a a'_{\al})a_{\al}$\ 
for $f\in \Hom_{-k}(A^{!}_{i+1}, A)$, $a\in A^{!}_i$.

Now, the entries of the Koszul complex are $A$-$k$-bimodules.
Let $k=\prod_y{\C_y}$ be a product of fields $\C$ and let
$y$ be the identities in these fields (local idempotents).
Let us fix $e\in k$ to be such a local idempotent. 
We can multiply the Koszul complex by $e$ from the right.
\begin{equation}\label{eKoszulComplex}
\cdots A\otimes{}^{\ast}(A^{!}_2)e\to  
A\otimes{}^{\ast}(A^{!}_1)e\to Ae\to \C_e,
\end{equation}
where $\C_e$ is the field corresponding to $e$ considered
as a simple left $A$-module. This is a projective resolution 
of $\C_e$.

Let $M$ be a finitely generated left module over $A$. Applying 
$\Hom_{A}(\underline{\ \ },M)$
to the complex (\ref{eKoszulComplex}) we get
the complex
\begin{equation}\label{eMKoszulComplex} 
\cdots e(A^{!}_2)\otimes M\leftarrow  e(A^{!}_1)\otimes M 
\leftarrow eM 
\end{equation}
calculating $\Ext^*_A(\C_e,M)$ with the differential
$d: e(A^{!}_i)\otimes M\to e(A^{!}_{i+1})\otimes M$
given by $d(ea\otimes m)=\sum ea a'_{\al} \otimes a_{\al}m$
for $a\in A^{!}_i$ and $m\in M$.

For a Koszul ring we have an algebra isomorphism 
\cite[Theorem 2.10.1]{BGS}
$$
\Ext^{*}_A(k,k)=(A^{!})^{\opp}
$$
which induces an isomorphism of $k$-$k$-bimodules
(with switched left and right action)
$$
\Ext^{i}_A(k,k)=A^{!}_i,
$$
and the isomorphism 
$$
\oplus_{y}\Ext^i_A(\C_e,\C_y)=
\Ext^{i}_A(\C_e,k)=eA^{!}_i,
$$
where $k=\prod_y{\C_y}$.

Rewriting the complex (\ref{eMKoszulComplex}) we get
\begin{equation}\label{ExteMKoszulComplex} 
\cdots \oplus_{y}
(\Ext^{2}_A(\C_e,\C_y)\otimes yM)\leftarrow  
\oplus_{y}(\Ext^{1}_A(\C_e,\C_y)\otimes yM) 
\leftarrow eM.
\end{equation} 

\sus{Flag varieties}\label{FlagKoszulAlg} 
Now let us apply this general machinery to our situation.
Recall that perverse sheaves on $X=G/P$ are modules over a Koszul 
algebra $A$ and $k=\prod_{W_{\SS}} \C_y$ where $W_{\SS}$ 
is the set of
Schubert strata, cf. \cite{BGS}. 
Let $e$ be the idempotent corresponding to the open 
Schubert stratum. Then, as a simple perverse sheaf $\C_e=\C_X[\dim X]$
is the shifted constant sheaf on the variety $X$. 

Let $M$ be a finitely generated left $A$-module and let us denote
by the same letter the corresponding perverse sheaf on $X$.
Let $\DD(X, W_{\SS})$ denote the derived category
of sheaves smooth along the Schubert stratification. 
We have 
\begin{equation}\label{EXThyperH}
\Ext^{i}_A(\C_e,M)= \Hom_{\DD(X, W_{\SS})}(\C_X[\dim X], M[i])
= \H^{i-\dim X}(M).
\end{equation}
Here the first equality is
due to
\cite[Corollary 3.3.2]{BGS} and the second is by definition
of hypercohomology. In particular,
$$
\Ext^{i}_A(\C_e,\C_y)= \H^{i-\dim X}(\ICS(\barr X_y))
=IH^{i-\dim X}(\barr X_y),
$$
where $X_y$ is the Schubert cell corresponding to the
idempotent $y\in k$.

Thus we can rewrite the complex (\ref{ExteMKoszulComplex})
as follows
\begin{equation}\label{IHeMKoszulComplex} 
\cdots 
\bigoplus_{y}  (IH^{1-\dim X}(\barr X_y) \otimes_k yM )
\leftarrow  \bigoplus_{y}(IH^{0-\dim X}(\barr X_y) \otimes_k yM )
\end{equation} 
with the differential 
$$
d: \bigoplus_{y}  (IH^{i}(\barr X_y) \otimes_k yM) 
\to \bigoplus_{y}(IH^{i+1}(\barr X_y) \otimes_k yM) 
$$
obtained as follows: 
\begin{equation}\label{differentialCoordinates}
d (c_y\otimes ym)= (\sum a'_{\al}\cdot c_y\otimes a_{\al} ym)
\end{equation}
for $c_y \in IH^{i}(\barr X_y)$, $m\in M$ and
$a'_{\al}\cdot c_y$ is the \emph{left} action of $\Ext^{*}_{A}(k,k)$
(identified with the \emph{right} action of $A^{!}$)
on $\oplus_{y}  IH^{i}(\barr X_y)$. Let us also
write down the formula for $d^2$:\

\begin{equation}\label{key}
d^2 (c\otimes m)=
\sum_{\al, \beta} (a'_{\beta}a'_{\al}\cdot c  
\otimes a_{\beta}a_{\al}\cdot m)
\end{equation}

\sus{Proposition}\label{MainProposition} 
The action of the tensor algebra 
$T(A_1) = k\oplus \bigoplus_{i=1}^{\infty} A_1^{\otimes i}$
on its finitely generated 
left module $M$ 
descends to the action of $A$ on $M$
if and only if the sequence of maps (\ref{IHeMKoszulComplex})
is a complex i.e., $d^2=0$.

First we prove the following two lemmas. 

\sus{Lemma}\label{Lemma1} 
Let $M$ be a finitely generated 
left $T(A_1)$-module. We have $d^2 =0$ for $d$   
constructed as in (\ref{IHeMKoszulComplex})
if and only if $(d')^2 = 0$ in the sequence of maps
\begin{equation}\label{ExtMKoszulComplex} 
\Ext^{2}_A(k,k)\otimes M\overset{d'}\leftarrow  
\Ext^{1}_A(k,k)\otimes M\overset{d'}  
\leftarrow M,
\end{equation} 
where 
$$
d'(m) = \sum_{\al} a'_{\al} \otimes a_{\al} m
\qquad\text{ and }
\qquad
d'(g\otimes m) = \sum_{\al} a'_{\al}\cdot g\otimes a_{\al} m
$$
for $g\in \Ext^{1}_A(k,k)$ and $m\in M$.

\begin{proof}
Indeed, 
\begin{equation}\label{keyCompare}
(d')^2 (m) = \sum_{\al, \beta} a'_{\beta}a'_{\al}
\otimes a_{\beta}a_{\al}\cdot m \in \Ext^{2}_A(k,k)\otimes M
\end{equation}

Consider a tensor product of the isomorphism of
the Theorem \ref{SoergelSatz} and the identity: 
$$
\phi: \Ext^{2}_A(k,k)\otimes_k M \to 
\Hom^2_{H^*(X)} (\oplus_y IH^*(\barr X_y),\oplus_y IH^*(\barr X_y))
\otimes_k M . 
$$
For $c\in \oplus_y IH^*(\barr X_y)$ let 
$$
\ev_c:  \Hom^2_{H^*(X)} (\oplus_y IH^*(\barr X_y),\oplus_y IH^*(\barr X_y))
\otimes_k M \to \oplus_y IH^*(\barr X_y)\otimes M
$$
be the map evaluating the first tensor multiple (i.e. $\Hom^2$) at $c$.
Comparing the formulas (\ref{key}) and (\ref{keyCompare})
we see that
$$
\ev_c(\phi( (d')^2 (m) )) = d^2(c\otimes m), 
$$
where $c\in \oplus_y IH^*(\barr X_y)$, $m\in M$, and
$d^2$ is constructed as in (\ref{key}). 
Since $\phi$ is an isomorphism, the statement of the lemma
follows.
\end{proof}

Recall that 
$A = T(A_1)/\left < R \right >$ where
$\left < R \right >$ is the ideal of relations
generated by $R\subset (A_1\otimes_k A_1)$.

\sus{Lemma}\label{Lemma2} 
The map $(d')^2$ defined in ($\ref{ExtMKoszulComplex}$) 
is equal to zero for a finitely generated left $T(A_1)$-module $M$ 
if and only if $r\cdot M =0$ for all $r\in R \subseteq (A_1\otimes A_1)$.

\begin{proof}  
Let us rewrite the complex 
(\ref{ExtMKoszulComplex}) as
\begin{equation}\label{RKoszulComplex} 
R^*\otimes M\overset{d'}\leftarrow  
A_1^!\otimes M\overset{d'}  
\leftarrow M,
\end{equation}
where $R^* = (A_1\otimes A_1)^*/R^{\perp} = A_2^{!}$.
We can give another description of the map $(d')^2$ as follows.
Denote $V = A_1 \otimes A_1$, and
$V^* = (A_1\otimes A_1)^*$. Consider the map
$$ 
\Id_{V^*} \otimes \operatorname{act}: V^*\otimes V\otimes M
\to V^* \otimes M,
$$
where $\operatorname{act}$ is the action of $V = A_1 \otimes A_1$
on $M$, and the tensor product of the quotient map 
$V^* \to V^*/R^{\perp}$ and $\Id_M$:
$$
p : V^*\otimes M \to V^*/R^{\perp} \otimes M = A_2^{!}\otimes M.
$$
Let $\psi = p\circ(\Id_{V^*} \otimes \operatorname{act})$.
Then for $m\in M$
$$
(d')^2 (m) = \psi
(\sum_{\al, \beta} (a'_{\al}\otimes a'_{\beta})
\otimes (a_{\beta}\otimes a_{\al}) \otimes m).
$$
(Remember that \emph{left} multiplication in $\Ext^{*}_A(k,k)$
is identified with the \emph{right} multiplication in $A^{!}$.)
Taking advantage of the identities $ V^*\otimes V = \Hom_k(V, V)$,
$V^* \otimes M = \Hom_k(V, M)$, and $A_2^{!}\otimes M = \Hom_k(R, M)$
we can consider $\psi$ as the composite map
$$
\Hom_k(V, V)\otimes M \to \Hom_k(V, M) \to \Hom_k(R, M).
$$
Now for $m\in M$
$$
(d')^2 (m) = \psi(\Id_V \otimes m) \in \Hom_k(R, M)
$$
and $(d')^2(m)(r) = r\cdot m$ for $r\in R$.
Thus $(d')^2 = 0$ if and only if $r\cdot m = 0$
for any $r\in R$ and any $m\in M$.
\end{proof}

\begin{proof}[Proof of Proposition \ref{MainProposition}] 
If the action of $T(A_1)$ descends
to the action of $A$, that is if $M$ is an $A$-module then
(\ref{IHeMKoszulComplex})
is a complex by construction. 

Suppose we are given a $T(A_1)$-module $M$
and suppose that (\ref{IHeMKoszulComplex})
is a complex. Then by Lemma \ref{Lemma1} and Lemma \ref{Lemma2}
we have $r\cdot m = 0$ for any $r\in R$ and any $m\in M$.
\end{proof}

\sus{From $T(A_1)$-modules to $\ic$-modules} Notice that
a finitely generated left 
(non-graded) $T(A_1)$-module $M$ for which the sequence of maps
(\ref{IHeMKoszulComplex}) is a complex
is exactly the same thing as a (non-mixed) $\ic$-module $V$
of finite length defined in \ref{mixedDefinition} and 
\ref{nonMixedDefinition}.
Indeed:

\sss{} The stalks $V_y = yM$ for all $y\in W_{\SS}$.
(Here we abuse notation by denoting a stratum of $X = G/P$
and the corresponding local idempotent by the same letter $y$.)
 
\sss{} The boundary maps $v(y, w)$, or rather their sum
$d = \oplus_{y, w} v(y, w)$ is given by
$$
\begin{aligned}
\bigoplus_y IH^{*}(X_y)\otimes yM & \simeq
\bigoplus_y IH^{*}(X_y)\otimes k \otimes yM \\ 
& \overset{\gamma_1}\to 
\bigoplus_y IH^{*}(X_y)\otimes (A_1^! \otimes A_1) \otimes yM \\
& \overset{\gamma_2}\to  \bigoplus_y IH^{*}(X_y)\otimes yM 
\end{aligned}
$$
where $\gamma_1$ is the canonical homomorphism 
$k \to A_1^! \otimes A_1$ tensored with $\Id's$;
and $\gamma_2$ is the tensor product of the \emph{right} action of
$A_1^!$ on $\oplus_y IH^{*}(X_y)$ (which is identified with
the \emph{left} action of $\Ext^1_A(k, k)$ on $\oplus_y IH^{*}(X_y)$
cf. (\ref{differentialCoordinates}))
and the action of
$A_1$ on $\oplus_y yM = M$.
\sss{} The Chain Complex Axiom is clearly equivalent to
the requirement that (\ref{IHeMKoszulComplex}) is a complex. 
In fact, (\ref{IHeMKoszulComplex}) is the total complex
of the $\ic$-module.

\sus{Last steps of the proof of The $d^2=0$ Theorem}
The complex (\ref{IHeMKoszulComplex}) calculates the hypercohomology 
of $M$ considered as a perverse sheaf due to the fact that 
the complex (\ref{eKoszulComplex}) is a projective resolution
of $\C_e$ and the formula \ref{EXThyperH}. The $d^2=0$ Theorem
is proved.

\sus{Example: $\P^1$ \cite{MV}} 
Let $X=\P^1$ stratified by two strata 
$X=S_0\sqcup S_1$ where $S_0=\C^0=\pt$, and $S_1=\C^1$.
We have
$$
IH^*(\barr S_0)=IH^*(\pt)=\C[0],
$$
where $\C[0]$ is the complex with zero differential 
and one-dimensional space
in position $0$, and
$$
IH^*(\barr S_1)=IH^*(\P^1)=\C[1]\oplus \C[-1],
$$
where $\C[1]\oplus \C[-1]$ is the complex with zero differential 
and two one-dimensional spaces
in positions $-1$ and $1$.
Now the spaces 
$$
\Hom^1_{H^*(\P^1)}(IH^*(\P^1),IH^*(\pt))
\ \text{ and }\  
\Hom^1_{H^*(\P^1)}(IH^*(\pt),IH^*(\P^1))
$$
are one dimensional. We will schematically depict a 
non-zero element
of the these spaces as
$$
\begin{CD}
        &             & \C[0] &        &  \\
          & \nearrow  &         &  &  \\
\C[1]& &&    &  \C[-1]  \\
\end{CD}
\qquad \text{   and   } \qquad
\begin{CD}
        &             & \C[0] &        &  \\
          &   &         &\searrow  &  \\
\C[1]& &&    &  \C[-1]  \\
\end{CD}
$$
respectively. Now we see that an $\ic$-module may be depicted as 
the following diagram 
$$
\begin{CD}
 &       &             & V_0 &        &  \\
 &         & a \nearrow   &         &  \searrow b &  \\
 &   V_1 &  &&    &  V_1  \\
\text{degree} :\ \ & -1 && 0 && 1 \\
\end{CD} \qquad\phantom{AAAA}
$$
with the relation $d^2=0$\  i.e.,\  $b\circ a=0$. 
Thus the category of $\ic$-modules
in this case is isomorphic to the category of quiver representations
\begin{equation}\label{qu}\nonumber
\qquad   \xymatrix{
        V_0 \ar@/_/[r]_{b} &  
        V_1  \ar@/_/[l]_{a} & 
}
\end{equation} 
with the relation $b\circ a=0$.

\se{Some combinatorics of Schubert varieties}\label{SchubertCombinatorics}

We start with some well known
combinatorial preliminaries. In this section 
$X=G/B$ is a full flag variety associated
to a simple algebraic group $G$ over $\C$. 

\sus{Cohomology ring}\label{scomb}
Let us denote $C=H^*(X)$. Following \cite{So98},
let $S(\fh)$ 
be the symmetric algebra over the Cartan $\fh$ of $\fg = \Lie G$,
let $S_{+}(\fh)\subset S(\fh)$ be the elements of positive 
degree and let $S_{+}(\fh)^W$ be the ideal generated by the 
$W$-invariants of $S_{+}(\fh)$. (Here the invariants are with respect to 
the usual ``linear'' action, cf. \cite{So98}.)
Recall that
$$
C = H^{*}(G/B) = S(\fh)/S_{+}(\fh)^W.
$$
Let $B^{-}$ be the opposite Borel. Denote by
$$
Y_w=B^{-}wB/B
$$
the opposite Schubert variety associated to $w\in W$,
and denote by $\si_w$ the class $[Y_w]$ in $H^{2l(w)}(X)$.
It is well known that the ring $H^{*}(X)$
is generated by the classes
$$
\si_{s_i},\qquad i=1, \dots, r,
$$
where $s_i$ are reflections
with respect to simple roots and $r$ is the rank of 
$\fg=\Lie G$. Moreover, we have the following
Chevalley's formula, cf. \cite{Ch}, for multiplying
the Schubert classes in $C = H^{*}(X)$
\begin{equation}\label{chevalley}
\si_{s_i}\cdot \si_w=
\sum_{\substack{\al \in \Phi^+\\ l(ws_{\alpha})=l(w)+1}} 
\langle\check\om_i, \al\rangle
\frac{(\al_i, \al_i)}{(\al, \al)} \si_{ws_{\al}},
\end{equation}
where $\check\om_i$ is the $i$-th 
fundamental coweight, and $\al_i$
is the $i$-th simple root, and 
$\langle\check\om_i, \al_j\rangle=\delta_{ij}$. 

\sus{}\label{subalgebras} 
Let $s_i$, $i\in \{1, \dots, r\}$ be a simple 
reflection in $W$ and let $C^{s_i}$ be the algebra
of invariants of $s_i$. It is a subalgebra of $C$.
We have the following description, explained to us 
by A. Postnikov, cf. \cite{Ch}:
$$
C^{s_i}=\spann\{ \si_w | ws_i > w,\  w\in W \}.
$$
Moreover, for any $s_i$, $i\in \{1, \dots, r\}$ we have
\begin{equation}\label{csdecomposition}
C= \si_{s_i}\cdot C^{s_i} \oplus C^{s_i}
\end{equation}

\sus{}\label{SoergelIndDesc}
 Now let us denote 
$$
\V_w=IH^*(\barr X_w).
$$
It is a graded $C$-module.
Due to Soergel, we have the following inductive 
description of $\V_w$. Take a reduced decomposition 
$w=s_1\dots s_l$
Denote
$$
C_{s_1\dots s_l}=C\otimes_{C^{s_l}} \otimes \dots C 
\otimes_{C^{s_1}} \C. 
$$
The action of $C$ of the leftmost tensor factor $C$
equips $C_{s_1\dots s_l}$ with a $C$-module
structure. By \cite{So90, So98} we know that
\begin{equation}\label{SoergelDecomposition}
C_{s_1\dots s_l}= \V_w \oplus \bigoplus_{y < w} \V^{n(y)}_y
\end{equation}
as a $C$-module.

\sss*{{\bf Example}} Type $A_2$.
\begin{enumerate}
\item $w=s_1s_2$:\ \ 
$
C_{s_1s_2}=\V_{w},
$
\item $w=s_1s_2s_1$:\ \
$
C_{s_1s_2s_1}=\V_{w}\oplus\V_{s_1},
\text{ and }
C_{s_2s_1s_2}=\V_{w}\oplus\V_{s_2}.
$
\end{enumerate}

\sus{} The $C$-module $C_{s_1\dots s_l}$ has a basis
\begin{equation}\label{tensorbasis}
a_l\otimes \dots \otimes a_1\otimes 1,
\end{equation}
where each
$$
a_i=\begin{cases}
\text{either} & 1, \\
\text{or} & \si_{s_i}. \\
\end{cases}
$$
In particular, 
$
\dim C_{s_1\dots s_l}=2^l.
$
This basis is homogeneous:
$$
\deg a_l\otimes \dots \otimes a_1\otimes 1=\deg a_l+\dots + \deg a_1,
$$
where $\deg 1=-1$ and $\deg \si_{s_i}=1$.
We know the action of 
$$
\si_{s_i},\ i=\{1, \dots, r \}
$$
on this basis by induction,
the decomposition \ref{csdecomposition},
and the Chevalley's 
formula (\ref{chevalley}).

\sus{} We can now give a description 
of $\ic$-modules purely in terms of combinatorics 
of the ring $C=H^*(X)$.
Recall that for each $w\in W$  we have 
a {\em graded} $C$-module
$\V_w=IH^*(\barr X_w)$ whose combinatorial description
is given in \ref{SoergelIndDesc}. Now a stalk of 
an $\ic$-module $M$ on $X=G/B$ is a vector space $M_w$
for every $w\in W$, and 
for every pair of $y, w\in W$ we have boundary maps
$$
m(y,w)\in
\Hom^1_{C}(\V_y,\V_w)\\
 \otimes \Hom(M_y, M_w)
$$
subject to the usual Chain Complex Axiom
(\ref{ChainComplexAxiom}). Let us 
denote this category of $\ic$-modules by 
$\AA_{\rm{Schubert}}(G/B)$.

\sus{Localization}\label{bbloc}
By the Beilinson-Bernstein Localization Theorem
(used in \cite{BB} to 
prove the Kazhdan-Lusztig conjecture \cite{KL79}, 
cf. \cite{KL80, BK})
and The $d^2=0$ Theorem
we obtain
\begin{theorem}\label{ictheorem}
The categories
$$
\OO_0\simeq \PP_{\rm{Schubert}}(G/B)\simeq \AA_{\rm{Schubert}}(G/B)
$$
are equivalent. 
\end{theorem}

\se{Quiver algorithm}\label{Algm}

Let us keep the setup of the previous section.
Using Theorem \ref{bbloc} 
we can calculate the relations between the arrows of 
the quiver for the category $\OO_0$ in the following way.

\sus{}\label{SectionQuiverd2}
 For $y, w \in W$ we have exactly $\mu(y, w)$ 
(note that $\mu(y, w)$ could be $0$)
arrows from $y$ to
$w$. Let us denote them by $a^1_{y, w}, \dots, a^{\mu(y, w)}_{y, w}$,
and consider these arrows as a basis in $A_1$, cf. Section 
\ref{KoszulComplexSection}.
Suppose we have a basis
in the vector space
$$
\Hom^1_{H^{*}(X)}(IH^{*}(\barr X_y),IH^{*}(\barr X_w))
=\Hom^1_C(\V_y, \V_w)
$$
given by matrices 
\begin{equation}\label{thematrices}
A^1_{y, w}, \dots, A^{\mu(y, w)}_{y, w}\in \Hom^1_C(\V_y, \V_w)
\end{equation}
with respect to some homogeneous bases in $\V_y$ and $\V_w$. 
Consider a matrix with elements in $T(A_1)$, 
cf. \ref{MainProposition}
$$
\widetilde d_{y, w} = \sum_{\substack{
i\in \{1, \dots, \mu(y, w)\}}} 
a^i_{y, w} A^i_{y, w}.
$$
Denote $\widetilde d = \sum_{y, w} \widetilde d_{y, w}$.
Then all the relations between the arrows are
encoded by the equation $\widetilde d^2=0$, 
or more precisely the 
linear relations between paths of length $2$  
going from $y$ to $w$
are the matrix elements of the matrix
$$
\sum_{\substack{
z\in W  \\
i\in \{1, \dots, \mu(z, w)\} \\
j\in \{1, \dots, \mu(y, z)\}}}
a^i_{z, w} A^i_{z, w}\circ a^j_{y, z} A^j_{y, z}=0,
$$
where $\circ$ is the matrix multiplication.
 
\sus{} Thus, all we need to get the relations
are the matrices (\ref{thematrices}) spanning
$$
\Hom^1_C(\V_y, \V_w).
$$
Algorithmically, it is enough to find homogeneous bases in 
$\V_w$ for $w\in W$ with the explicit action of 
the generators 
$$
\si_{s_i}, \qquad i=\{1\dots r\}
$$
of $C$. 
Then we can take all linear maps of degree $1$ between 
$\V_y$ and $\V_w$ and solve a system of linear equations
to find those linear maps which commute with the action 
of (the generators of) $C$.

\sus{}\label{HomBasisVw} 
We will find a homogeneous basis in $\V_w$ by induction on
the length of the element $w\in W$.

\sss*{{\bf Step $0$}} First of all, if $w=1$, then $\V_1=\C$.

\sss*{{\bf Step $1$}} By induction, suppose that we already 
have homogeneous bases with the explicit action of 
(the generators of) 
$C$ in
$\V_y$ for $y <  w$. Take a reduced decomposition 
$w=s_1\dots s_l$. It is easy to see from Theorem
\ref{SoergelSatz} that
$$
\Hom^0_C(\V_y, \V_w)=
\begin{cases}
\C, & y=w \\
0, & \text{otherwise}
\end{cases}
$$
From that and the decomposition
\ref{SoergelDecomposition}
we can find matrices 
(with respect to the induction-assumed basis
in $\V_y$ and the basis (\ref{tensorbasis})
in $C_{s_1\dots s_l}$)
spanning the space
$$
\Hom^0_C(\V_y, C_{s_1\dots s_l})
$$
by considering degree $0$ maps between
$\V_y$ and $C_{s_1\dots s_l}$ and making sure they
commute with the action of the generators of $C$
by solving a system of linear equations. 
Applying these matrices to the induction-assumed basis
in $\V_y$ we get a 
homogeneous linearly independent system of vectors
$
v_1, \dots, v_m
$
in 
$C_{s_1\dots s_l}$ spanning the subspace
$$
Y=\bigoplus_{y < w} \V^{n(y)}_y.
$$

\sss*{{\bf Step $2$}} 
Now we have a $C$-module $U=C_{s_1\dots s_l}$ with the basis
(\ref{tensorbasis}) which we will number in an 
arbitrary way as
$x_1, \dots, x_{2^l}$, and its $C$-submodule 
$$
Y=\bigoplus_{y < w} \V^{n(y)}_y
$$
with the basis 
$v_1, \dots, v_m$ which we know in terms
of the basis $x_1, \dots, x_{2^l}$:
$$
v_i=\sum_{k} \al^k_i x_k.
$$

We need to find a basis in the $C$-module $U/Y$
with the explicit action of (the generators of) $C$. 
This is a problem of computational commutative algebra
algorithmic solutions to which are well known, cf. e.g. 
\cite[Chapter 15]{E}. Let us just outline the algorithm:
\begin{itemize}
\item take a basis element $x_i$ not in $Y$ and generate
a $C$-submodule $Cx_i \subset U$; 
\item record a basis in $Cx_i - Y$
consisting of (linear combinations of) the
elements of the form $\si_w x_i$, $w\in W$
\item if $Cx_i + Y \neq U$,
take another basis element $x_j$ not in $Cx_i + Y$,
and generate a $C$-submodule $Cx_j$.
\item proceed until  $Cx_i + Y = U$. The algorithm will stop
since $\dim_{\C} U = 2^l < \infty$.
\end{itemize}
This algorithm involves arbitrary choices of
basis elements $x_i, x_j, \dots$.

\sss*{{\bf Step $3$}} 
Proceeding by induction we find a
homogeneous basis in $\V_w$, for all $w\in W$.
Due to the choices above this basis will depend on the
particular realization of the algorithm.

\sus{Remark} In order to improve the performance of an actual
computer realization of this algorithm, one should use
a number of shortcuts. For example, instead of using 
the module $C_{s_1\dots s_l}$ to extract a basis in $\V_w$,
we could use $C\otimes_{C^s}\V_{w'}$, where $w=w's$ and 
$l(w)=l(w')+1$. (I am grateful to 
T.~ Braden for this and other suggestions.)  
A detailed implementation of the algorithm in the $A_2$ case is worked 
out in the Appendix (Section 8).

\sus{Remark: one quiver = all quivers} Recall the Koszul algebra $A$ 
underlying the category $\OO_0\simeq\PP_{\rm{Schubert}}(G/B)$
from \cite{BGS} and \ref{FlagKoszulAlg}.
The algebra
$A$ is generated over $A_0$ by $A_1$ and it is quadratic, that is
the ideal of relations is generated by an $A_0\text{-}A_0$-bimodule
$$
R\subset A_1\otimes_{A_0} A_1.
$$
The quiver of the category $\OO_0$ in this language is just 
a choice of basis 
in $A_1$ (arrows) and a choice of basis in the subspace $R$
(relations between arrows). Once we have a 
basis in $A_1$ and a basis in $R$ 
we can obtain any other bases in $A_1$ 
and $R$ by linear transformations.

\se{Weight filtration and the functor $\HH$}\label{WeightSec}

In this section we construct a functor $\HH$ from
mixed perverse sheaves to mixed $\ic$-modules.

\sus{} Let $\PP^{\mixed}_{\SS}(X)$ be a mixed category of 
mixed perverse sheaves. Examples of such categories include
\cite[4.4, 4.5]{BGS} and \cite{V00, Br03}.
Let $v: \PP^{\mixed}_{\SS}(X)\to \PP_{\SS}(X)$
be the degrading functor.

By \cite[Lemma 4.1.2]{BGS}, cf. \cite{BBD}, every object 
$A\in \PP^{\mixed}_{\SS}(X)$ has a unique finite decreasing 
weight filtration $W_{\bullet}=W_{\bullet}A$
such that $W_mA/W_{m-1}A$ is
a pure object of weight $m$.

\sss{} Suppose that we have a functor 
$K: \PP^{\mixed}_{\SS}(X) \to C^b(\Vect)$
to the category of bounded complexes of (graded)
$R$-modules such that for $A\in \PP^{\mixed}_{\SS}(X)$
the cohomology of $K(A)$ is the hypercohomology of $A$
and the weight filtration of $A$ induces a decreasing filtration
\begin{equation}\label{wfilt}
K(A)=K(A)\supseteq \cdots \supseteq K(W_mA)
\supseteq  K(W_{m-1}A) \cdots 
\end{equation}
on $K(A)$, and moreover, the cohomology of $K(W_mA)/K(W_{m-1}A)$
is the hypercohomology of $W_mA/W_{m-1}A$.

\sus{} Assuming the statement of The $d^2=0$ Theorem, one could take 
the total complex of the corresponding (non-mixed)
$\ic$-module corresponding to $v(A)$ as the complex $K(A)$ for 
$A\in \PP^{\mixed}_{\SS}(X)$. 
In particular, we can do it
for simplicial complexes and flag varieties. We will
fix this choice of $K$ for the rest of this section.

Consider now the first term $E^1=E^1(A)$ of the spectral sequence
associated to the filtration (\ref{wfilt}). One could look at $E^1$
as a mixed $\ic$-module. Indeed, for $A\in \PP^{\mixed}_{\SS}(X)$
the stalks of the corresponding $\ic$-module $V=\HH(A)$
are determined from the formula
$$
H^i \left(K(W_mA)/K(W_{m-1}A)\right)=\bigoplus_S 
V^{-m}_S\otimes IH^i(\barr S)
$$
and the boundary maps of $V$ are the differentials of $E^1$.
The total complex of $V$ constructed in this way
is precisely the diagonal complex of $E^1$.

Thus we have constructed a functor 
$\HH:\PP^{\mixed}_{\SS}(X)\to\AA^{\mixed}_{\SS}(X)$.
Observe that the spectral sequence degenerates 
at the second term and so the total complex 
of the $\ic$-module $\HH(A)$ calculates the hypercohomology of $A$
for $A\in \PP^{\mixed}_{\SS}(X)$. In other words we have the following

\begin{theorem} Assuming the statement of {\rm{The $d^2=0$ Theorem}},
the total complex of the mixed $\ic$-module
calculates the hypercohomology of the
corresponding mixed perverse sheaf.
\end{theorem}

\sus{} 
Let $\widetilde\PP(G/B)$
be the mixed category of mixed perverse sheaves considered in
\cite[4.4]{BGS}. Let $\widetilde\AA(G/B)$ 
be the subcategory of the category $\AA^{\mixed}(G/B)$
with coefficients in $\Q_l$
consisting of objects $V$ such that 
$$
V^i_S=0 \qquad \text{ unless }\qquad \dim S=i\mod 2.
$$
Then the functor 
$
\HH: \widetilde\PP(G/B)\to\widetilde\AA(G/B)
$ 
provides the mixed version of the equivalence established
in Section \ref{SecFlagProof}.

\se{Further directions and final remarks}\label{FinalRemarks}

\sus{A combinatorial challenge} The basis in $\V_w$ we constructed 
in \ref{HomBasisVw} is very 
non-canonical: it depends on the reduced decompositions of $w$ 
we choose for each $w$, and on the way we number our bases. 
However, once we have found \emph{one} basis in $\V_w$ we can obtain 
\emph{any other} basis
by a (grading-preserving) linear transformation. 

It would be very interesting to construct a
\emph{distinguished} homogeneous basis in the $C$-module $\V_w$
for all $w\in W$ independent of any choices. The basis
in question should specialize to the basis of Schubert
cycles when the variety $\barr X_w$ is smooth, and
the formulas for the action of $C$
should generalize the Chevalley's formula (\ref{chevalley}).

\sus{Moment graphs} In \cite{BM} T.~ Braden and R. ~MacPherson
give another construction of a module structure on $IH^{*}(\barr X_w)$ 
and $IH^{*}_T(\barr X_w)$, where the latter is the 
$T$-equivariant intersection homology. One could use their techniques
to compute with $\ic$-modules. It would also be very interesting
to define \emph{$T$-equivariant} $\ic$-modules,
perhaps related to $T$-equivariant perverse sheaves
and \emph{singular} blocks of $\OO$.

\sus{Parabolic and singular case} In Section \ref{Algm} 
we have considered 
only the full flag variety
and the regular block of the category $\OO_0$.
However, the algorithm goes through without change
for a parabolic flag variety $G/P$
and the corresponding regular parabolic block of the 
category $\OO$. The parabolic Chevalley's formula is known, 
cf. \cite{FW}. The quiver for a singular block 
can be obtained by Koszul (and thus quadratic) duality
from the corresponding parabolic quiver, cf \cite{BGS}.

\sus{Harish-Chandra modules} In fact, in 
\cite[Section 1]{Gelf} \linebreak
I.~ M.~ Gelfand asks for a quiver describing
(a block of) the category of Harish-Chandra modules
over an arbitrary simple group. Harish-Chandra modules 
reduce to category $\OO$ in some cases \cite{BernGelf}.
We hope that our 
methods will provide an answer to that question 
as well, and help prove Soergel's Langlands duality 
conjecture, cf. \cite{So01}. (There are even more general
representation categories, related to arbitrary Coxeter systems,
for which the theory partially goes through, cf. \cite{D}.)

\sus{\bf Toric varieties}\label{toric}
First we introduce the \emph{relative} $\ic$-modules. 
Let $D(X)=D^b_c(X)$ be 
the bounded derived category of constructible sheaves on $X$.
We could define a relative version of $\ic$-modules
by using another $\partial$-functor 
$D^b_c(X)\to D^b(A)$, where $A$ is an abelian category,
instead of hypercohomology. A good example of such a functor would be
the derived functor $Rf_*: D^b_c(X)\to  D^b_c(Y)$ 
associated to a map $f:X\to Y$.

It would be interesting to
look at
the relative $\ic$-modules related to
the Koszul category
of perverse sheaves on toric varieties studied in \cite{Br03}. 
Let $X$ be a toric variety
with the action of the torus $T$ and let $\mu: X\to X/T$
be the quotient map. The relative $\ic$-modules with respect to
$\mu$ would still have vector spaces $V_{\sigma}$ 
as stalks at each face $\sigma$,
and
elements of the space
$$
\Hom^1_{\mu_*R_X}(\mu_*\ICS({\sigma}, V_{\sigma}), 
\mu_*\ICS({\tau}, V_{\tau}))
$$
for any two faces ($T$-orbits) $\sigma$ and $\tau$
as boundary maps, satisfying the $d^2=0$ axiom.
(Here
$\mu_*$ is the functor between derived categories, and $R_X$
is the constant sheaf on $X$.) The objects 
$\mu_*R_X$, $\mu_*\ICS_{\sigma}$, and $\mu_*\ICS_{\tau}$ 
have a very nice combinatorial description, cf. \cite[1.3]{Br03}. 

\sus{\bf Affine Grassmannians} The $d^2=0$ Theorem trivially holds
for the semisimple category of $G(O)$-equivariant
perverse sheaves of $\C$-vector spaces on the affine Grassmannian
$G(K)/G(O)$. It would be interesting to prove an analogous result
for a suitable Koszul category of perverse sheaves
of vector spaces over a field of positive characteristic.

\sus{\bf Saper's $L$-modules}\label{SapSec}
In a remarkable recent paper \cite{S} L.~ Saper
introduced the notion of $L$-modules, in order
to prove a conjecture of Rapoport and Goresky-MacPherson.
It seems that $L$-modules should be related to 
$\ic$-modules on the (compactifications of) locally symmetric spaces.
It would be very interesting to understand the relationship between
these structures and perverse sheaves.

\sus{\bf Generalized $\OPH$-modules} 
If we take the usual homology instead of intersection homology 
in the definition of $\ic$-modules, 
we will get the usual 
triangulation-constructible sheaves on simplicial complexes,
cf. \cite{KS, M93, V00}. Thus we have a machine transforming
a homology theory into an abelian 
(or, more generally, $A_{\infty}$) category of ``sheaves.''
It would be
intriguing to consider $\ic$-modules with an arbitrary 
(generalized) homology theory $\OPH$ instead of 
intersection homology.

\se{Appendix. The quiver in the $A_2$ case}

The quiver relations in the $A_2$ case were obtained in 
\cite{I82, St03}, but we would still like to illustrate
our quiver algorithm from section \ref{Algm} 
on this example. The quiver we obtain 
is defined ``over the integers'', cf. 
\ref{openQuestionPMone}. 

\sus{Cohomology ring} In the $A_2$ case the Weyl group
$W = {\mathfrak S_3} = \{1, s_1, s_2, s_1s_2, s_2s_1, s_1s_2s_1\}$ 
and the cohomology ring $C = H^{*}(G/B)$ has a basis
of $6$ elements: $C = \spann\{1, \si_{s_1}, \si_{s_2}, 
\si_{s_1s_2}, \si_{s_2s_1}, \si_{s_1s_2s_1} \}$.
For simplicity we denote $\si_i = \si_{s_i}$, $i=1, 2$.

Using the Chevalley's formula (\ref{chevalley})
we compute the following partial multiplication table
for this algebra:
\vskip .2in
\begin{center}
\begin{tabular}{|l|c|c|}
\hline
 & $\si_1$ & $\si_2$ \\
\hline
$1$ & $\si_1$ & $\si_2$ \\
\hline
$\si_1$ & $\si_{s_2s_1}$ & $\si_{s_2s_1} + \si_{s_1s_2}$ \\
\hline
$\si_2$ & $\si_{s_2s_1} + \si_{s_1s_2}$ & $\si_{s_1s_2}$ \\
\hline
$\si_{s_1s_2}$ & $\si_{s_1s_2s_1}$ & $0$  \\
\hline
$\si_{s_2s_1}$ & $0$ & $\si_{s_1s_2s_1}$ \\
\hline
$\si_{s_1s_2s_1}$ & $0$ & $0$  \\
\hline
\end{tabular}
\end{center}
\vskip .2in

\sss{} The subalgebras $C^{s_i}$, $i=1, 2$ presented in 
subsection \ref{subalgebras} can be explicitly described as follows:
$$
C^{s_1} = \spann\{1, \si_{2}, \si_{s_1s_2} \} ,
$$
and
$$
C^{s_2} = \spann\{1, \si_{1}, \si_{s_2s_1}\} .
$$

\sss{} It is easy to see the decomposition (\ref{csdecomposition})
in this case:
$$
\si_1 \cdot C^{s_1} \oplus C^{s_1}
= \spann\{\si_{1}, \si_{s_2s_1} + \si_{s_1s_2}, \si_{s_1s_2s_1}\}
\oplus \spann\{1, \si_{2}, \si_{s_1s_2}\} = C ,
$$
and
$$
\si_2 \cdot C^{s_2} \oplus C^{s_2}
= \spann\{\si_{2},  \si_{s_2s_1} + \si_{s_1s_2},  \si_{s_1s_2s_1}\}
\oplus \spann\{1, \si_{1}, \si_{s_2s_1}\} = C .
$$

\sus{The modules $\V_w$}\label{modulesVw} 
These modules 
can be explicitly described 
as follows:

\sss{} $\V_{1} = \C$, and $\si_i$, $i=1, 2$ act by $0$.

\sss{} $\V_{s_1} = C\otimes _{C^{s_1}} \C 
=\spann\{1\otimes 1, \si_1\otimes 1\}$.
Action of $\si_1$:
$$
\si_1\cdot 1\otimes 1 = \si_1 \otimes 1
$$
and
$$
\si_1\cdot (\si_1\otimes 1) = \si_1^2 \otimes 1 =
\si_{s_2s_1} \otimes 1 = (\si_1\si_2 - \si_{s_1s_2})\otimes 1
=  \si_1\otimes \si_2\cdot 1 - 1\otimes \si_{s_1s_2}\cdot 1 = 0.
$$
Action of $\si_2$:
$$
\si_2\cdot 1\otimes 1 = \si_2 \otimes 1 = 1\otimes \si_2\cdot 1 = 0
$$
and
$$
\si_2\cdot (\si_1\otimes 1) = \si_1\si_2\otimes 1 = 
\si_1 \otimes \si_2\cdot 1 = 0.
$$
Thus the action of $\si_1$ and $\si_2$ in this basis is given by
the matrices
$$
\si_1 = 
\left(
\begin{matrix}
0 & 0 \\
1 & 0
\end{matrix} 
\right ) ,
\qquad 
\si_2 = 
\left(
\begin{matrix}
0 & 0 \\
0 & 0
\end{matrix}
\right ) .
$$

\sss{} $\V_{s_2} = C\otimes _{C^{s_2}} \C 
=\spann\{1\otimes 1, \si_2\otimes 1\}$.
The action of $\si_1$ and $\si_2$ in this basis is given by
the matrices
$$
\si_1 = 
\left(
\begin{matrix}
0 & 0 \\
0 & 0
\end{matrix}
\right ) ,
\qquad 
\si_2 = 
\left(
\begin{matrix}
0 & 0 \\
1 & 0
\end{matrix}
\right ) .
$$

\sss{} $\V_{s_1s_2} = C\otimes _{C^{s_2}} C\otimes _{C^{s_1}}\C$.
Note that 
$$
\V_{s_1s_2} =\spann\{1\otimes 1\otimes 1, \si_2\otimes 1\otimes 1,
1\otimes  \si_1\otimes 1,  \si_2\otimes \si_1\otimes 1\} .
$$
The action of $\si_1$ and $\si_2$ in this basis is given by the matrices:
$$
\si_1 = 
\left( 
\begin{matrix} 
0 & 0 & 0 & 0 \\ 
0 & 0 & 0 & 0 \\
1 & 0 & 0 & 0 \\
0 & 1 & 0 & 0 \\
\end{matrix}
\right ) ,
\qquad
\si_2 = 
\left( 
\begin{matrix} 
0 & 0 & 0 & 0 \\ 
1 & 0 & 0 & 0 \\
0 & 0 & 0 & 0 \\
0 & 1 & 1 & 0 \\
\end{matrix}
\right ) .
$$

\sss{} $\V_{s_2s_1} = C\otimes _{C^{s_1}} C\otimes _{C^{s_2}}\C$.
Note that 
$$
\V_{s_1s_2} =\spann\{1\otimes 1\otimes 1, \si_1\otimes 1\otimes 1,
1\otimes  \si_2\otimes 1,  \si_1\otimes \si_2\otimes 1\} .
$$
The action of $\si_1$ and $\si_2$ in this basis is given by the matrices:
$$
\si_1 = 
\left( 
\begin{matrix} 
0 & 0 & 0 & 0 \\ 
1 & 0 & 0 & 0 \\
0 & 0 & 0 & 0 \\
0 & 1 & 1 & 0 \\
\end{matrix}
\right ) ,
\qquad
\si_2 = 
\left( 
\begin{matrix} 
0 & 0 & 0 & 0 \\ 
0 & 0 & 0 & 0 \\
1 & 0 & 0 & 0 \\
0 & 1 & 0 & 0 \\
\end{matrix}
\right ) .
$$

\sss{} $C_{s_1s_2s_1} = 
C\otimes _{C^{s_1}} C\otimes _{C^{s_2}}C\otimes _{C^{s_1}}\C$
We have
$$
\begin{aligned}
C_{s_1s_2s_1}  = 
\spann\{
& 1\otimes 1\otimes 1\otimes 1, 
1\otimes 1\otimes \si_1\otimes 1, \\
& 1\otimes \si_2 \otimes 1\otimes 1,
\si_1\otimes 1\otimes 1\otimes 1, \\
& 1\otimes \si_2\otimes \si_1\otimes 1,
\si_1\otimes 1\otimes \si_1\otimes 1, \\
& \si_1\otimes \si_2\otimes 1\otimes 1,
\si_1\otimes \si_2\otimes \si_1\otimes 1
\}
\end{aligned}
$$

The action table is as follows:

\vskip .2in
\begin{center}
\begin{tabular}{|l|c|c|}
\hline
 & $\si_1$ & $\si_2$ \\
\hline
$1\otimes 1\otimes 1\otimes 1$ 
& $\si_1\otimes 1\otimes 1\otimes 1$  
& $1\otimes\si_2 \otimes 1\otimes 1$ \\
\hline
$1\otimes 1\otimes \si_1\otimes 1$
& $\si_1\otimes 1\otimes \si_1\otimes 1$
& $1\otimes \si_2\otimes \si_1\otimes 1$ \\
\hline
$1\otimes  \si_2\otimes 1\otimes 1$ 
& $\si_1\otimes \si_2\otimes 1\otimes 1$ 
& $1\otimes  \si_2\otimes \si_1\otimes 1$ \\
\hline
$\si_1\otimes 1\otimes 1\otimes 1$ 
&  $\si_1\otimes \si_2\otimes 1\otimes 1$ -  
$1\otimes  \si_2\otimes \si_1\otimes 1$
& $\si_1\otimes \si_2\otimes 1\otimes 1$ \\
\hline
$1\otimes  \si_2\otimes \si_1\otimes 1$
& $\si_1\otimes \si_2\otimes \si_1\otimes 1$
& 0 \\
\hline
$\si_1\otimes 1\otimes \si_1\otimes 1$ &
$\si_1\otimes \si_2\otimes \si_1\otimes 1$
& $\si_1\otimes \si_2\otimes \si_1\otimes 1$ \\
\hline
$\si_1\otimes \si_2\otimes 1\otimes 1$
& $\si_1\otimes \si_2\otimes \si_1\otimes 1$
& $\si_1\otimes \si_2\otimes \si_1\otimes 1$ \\
\hline
$\si_1\otimes \si_2\otimes \si_1\otimes 1$ & 0 & 0 \\
\hline
\end{tabular}
\end{center}
\vskip .2in

The decomposition (\ref{SoergelDecomposition}) in this case is
$C_{s_1s_2s_1} = \V_{s_1s_2s_1} \oplus \V_{s_1}$.
The submodule isomorphic to $\V_{s_1}$ is spanned in 
$C_{s_1s_2s_1}$ by:
$$
\begin{aligned}
\V_{s_1} & \hookrightarrow C_{s_1s_2s_1} \\ 
1\otimes 1 & \mapsto 1\otimes 1\otimes \si_1\otimes 1-
 1\otimes \si_2\otimes 1\otimes 1 \\
\si_1\otimes 1 & \mapsto \si_1\otimes 1\otimes \si_1\otimes 1-
\si_1\otimes \si_2\otimes 1\otimes 1
\end{aligned}
$$
Following Step 2 of section \ref{HomBasisVw}, choose
a basis element $1\otimes 1\otimes 1\otimes 1 \in C_{s_1s_2s_1}$
and generate a $C$-submodule 
$C\cdot (1\otimes 1\otimes 1\otimes 1) \subset C_{s_1s_2s_1}$.
We have
$$
\begin{aligned}
C\cdot (1\otimes 1\otimes 1\otimes 1) & = \spann\{
1\otimes 1\otimes 1\otimes 1, \\
& \si_1\otimes 1\otimes 1\otimes 1,\\
& 1\otimes \si_2\otimes 1\otimes 1,\\
& \si_1\otimes \si_2\otimes 1\otimes 1 -  
 1\otimes  \si_2\otimes \si_1\otimes 1, \\
& 1\otimes  \si_2\otimes \si_1\otimes 1 \\
& \si_1\otimes \si_2\otimes \si_1\otimes 1
\end{aligned}
$$
We see that $C\cdot (1\otimes 1\otimes 1\otimes 1) + \V_{s_1} =
C_{s_1s_2s_1}$,
so $\V_{s_1s_2s_1} =  C\cdot (1\otimes 1\otimes 1\otimes 1)$
with the basis as above. The matrices of 
$\si_1$ and $\si_2$ in this basis are as follows:
$$
\si_1 = 
\left( 
\begin{matrix} 
0 & 0 & 0 & 0 & 0 & 0 \\ 
1 & 0 & 0 & 0 & 0 & 0 \\
0 & 0 & 0 & 0 & 0 & 0 \\
0 & 0 & 1 & 0 & 0 & 0 \\
0 & 1 & 1 & 0 & 0 & 0 \\
0 & 0 & 0 & 1 & 0 & 0
\end{matrix}
\right ) ,
\qquad
\si_2 = 
\left( 
\begin{matrix} 
0 & 0 & 0 & 0 & 0 & 0 \\ 
0 & 0 & 0 & 0 & 0 & 0 \\
1 & 0 & 0 & 0 & 0 & 0 \\
0 & 1 & 1 & 0 & 0 & 0 \\
0 & 1 & 0 & 0 & 0 & 0 \\
0 & 0 & 0 & 0 & 1 & 0
\end{matrix}
\right ) .
$$

\sss{} Notice that for $\V_{s_1s_2s_1} = C$ as a $C$-module, 
so we could also use the basis of ``Schubert cycles'' 
$\{1, \si_{s_1}, \si_{s_2}, 
\si_{s_1s_2}, \si_{s_2s_1}, \si_{s_1s_2s_1} \}$. 

\sus{Matrices representing $\Hom^1$} 

\sss{Setup} For any pair $w, y \in W = {\mathfrak S}_3$.
We need to find matrices representing homomorphisms 
$f_{w, y}: \V_w \to \V_y$ of degree $1$
such that $\si_i f = f\si_i$ for $i=1, 2$. 
We use the bases in $\V_w$ constructed the section \ref{modulesVw}.

We will collect these matrices into one matrix
of the map $d: \oplus_w\V_w \to \oplus_w \V_w$.
Since $\dim \oplus_w\V_w = 19$, the matrix of $d$
has dimension $19\times 19$, and it has $16$ nonzero
blocks of the form $f_{w, y}$, $ w, y, \in {\mathfrak S}_3$.
\begin{equation}
d = \left(
\begin{matrix}
0 & f_{s_1, 1} & f_{s_2, 1} & 0 & 0 & 0 \\
f_{1, s_1} & 0 & 0 & f_{s_1s_2, s_1} & f_{s_2s_1, s_1} & 0 \\
f_{1, s_2} & 0 & 0 & f_{s_1s_2, s_2} & f_{s_2s_1, s_2} & 0 \\
0 & f_{s_1, s_1s_2} &  f_{s_2, s_1s_2} & 0 & 0 &  f_{s_1s_2s_1, s_1s_2} \\
0 & f_{s_1, s_2s_1} &  f_{s_2, s_2s_1} & 0 & 0 &  f_{s_1s_2s_1, s_2s_1} \\
0 & 0 & 0 & f_{s_1s_2, s_1s_2s_1} &  f_{s_2s_1, s_1s_2s_1} & 0
\end{matrix}
\right ) .
\end{equation}

\sss{The matrices  $f_{w, y}$, $ w, y, \in {\mathfrak S}_3$}
We record the $16$ matrices $f_{w, y}$ obtained by solving
the linear systems $\si_i f = f\si_i$ for $i=1, 2$.

$$
f_{1,s_1} =
\left (
\begin{matrix}
0 \\
1
\end{matrix}
\right )
\qquad
f_{1,s_2} =
\left (
\begin{matrix}
0 \\
1
\end{matrix}
\right )
\qquad
f_{s_1, 1} =
\left (
\begin{matrix}
1 & 0 
\end{matrix}
\right )
\qquad
f_{s_2, 1} =
\left (
\begin{matrix}
1 & 0
\end{matrix}
\right )
$$

$$
f_{s_1, s_1s_2}
=\left(
\begin{matrix}
0 & 0 \\
1 & 0 \\
-1 & 0 \\
0 & 1
\end{matrix}
\right)
\qquad
f_{s_1, s_2s_1}
=\left(
\begin{matrix}
0 & 0 \\
0 & 0 \\
1 & 0 \\
0 & 1 \\
\end{matrix}
\right)
$$

$$
f_{s_2, s_1s_2}
=\left(
\begin{matrix}
0 & 0 \\
0 & 0 \\
1 & 0 \\
0 & 1
\end{matrix}
\right)
\qquad
f_{s_2, s_2s_1}
=\left(
\begin{matrix}
0 & 0 \\
1 & 0 \\
-1 & 0 \\
0 & 1
\end{matrix}
\right)
$$

$$
f_{s_1s_2, s_1}
=\left(
\begin{matrix}
1 & 0 & 0 & 0 \\
0 & 0 & 1 & 0
\end{matrix}
\right)
\qquad
f_{s_2s_1, s_1}
=\left(
\begin{matrix}
1 & 0 & 0 & 0 \\
0 & 1 & 0 & 0
\end{matrix}
\right)
$$

$$
f_{s_1s_2, s_2}
=\left(
\begin{matrix}
1 & 0 & 0 & 0 \\
0 & 1 & 0 & 0
\end{matrix}
\right)
\qquad
f_{s_2s_1, s_1}
=\left(
\begin{matrix}
1 & 0 & 0 & 0 \\
0 & 0 & 1 & 0
\end{matrix}
\right)
$$

$$
f_{s_1s_2, s_1s_2s_1}
=\left(
\begin{matrix}
0 & 0 & 0 & 0 \\
1 & 0 & 0 & 0 \\
0 & 0 & 0 & 0 \\
0 & 1 & 0 & 0 \\
0 & 1 & 1 & 0 \\
0 & 0 & 0 & 1
\end{matrix}
\right)
\qquad
f_{s_2s_1, s_1s_2s_1}
=\left(
\begin{matrix}
0 & 0 & 0 & 0 \\
0 & 0 & 0 & 0 \\
1 & 0 & 0 & 0 \\
0 & 1 & 1 & 0 \\
0 & 1 & 0 & 0 \\
0 & 0 & 0 & 1
\end{matrix}
\right)
$$

$$
f_{s_1s_2s_1, s_1s_2}
=\left(
\begin{matrix}
1 & 0 & 0 & 0 & 0 & 0 \\
0 & 0 & 1 & 0 & 0 & 0 \\
0 & 1 & 0 & 0 & 0 & 0 \\
0 & 0 & 0 & 1 & 0 & 0 
\end{matrix}
\right)
\qquad
f_{s_1s_2s_1, s_2s_1}
=\left(
\begin{matrix}
1 & 0 & 0 & 0 & 0 & 0 \\
0 & 1 & 0 & 0 & 0 & 0 \\
0 & 0 & 1 & 0 & 0 & 0 \\
0 & 0 & 0 & 0 & 1 & 0 
\end{matrix}
\right)
$$

\sus{The quiver notation}

\sss{Vertex enumeration} The quiver of the algebra in question looks as 
follows (here we actually have two arrows in the opposite directions
between every two vertices connected by the graphical image $\leftrightarrow$)
 
\begin{equation}
\xymatrix{
& 6 \ar@{<->}[dr] \ar@{<->}[dl]& \\
4 \ar@{<->}[d]\ar@{<->}[drr]  & & 5 \ar@{<->}[d] \ar@{<->}[dll]\\
2 \ar@{<->}[dr]  & & 3 \ar@{<->}[dl]\\
& 1 &
}
\end{equation}
where in order to harmonize the notation with \cite{I82, St03}
we number the elements of $\mathfrak S_3$ as follows:
\begin{enumerate}
\item $s_1s_2s_1$ is number $1$, 
\item $s_1s_2$ is number $2$,
\item $s_2s_1$ is number $3$,
\item $s_1$ is number $4$,
\item $s_2$ is number $5$,
\item $1$ is number $6$.
\end{enumerate}

\sss{Arrow notation} 
We have $16$ arrows on the quiver. The arrows will be denoted
by specifying the numbers of their initial and terminal vertices.
For example the arrow from $s_1$ to $1$ (vertex number $4$ to
vertex number $6$) will be denoted as $(46)$. 

\sss{Path notation} A path
on the quiver will be denoted by the sequence of the
vertices through which it goes. For example the path
of length $2$ starting at $s_1s_2s_1$ (vertex number $1$) 
going through $s_1s_2$ (vertex number $2$) and ending at 
$s_1$ (vertex number $4$) will be denoted by the sequence $(124)$.

\sss{The free quiver algebra} Consider the algebra freely
generated by the arrows of the quiver above, i.e.
the tensor algebra 
$T(A_1) = k\oplus \bigoplus_{i=1}^{\infty} A_1^{\otimes i}$
where $k$ is the semisimple ring spanned by the 
idempotents at the vertices and
$A_1$ is the $k-k$-bimodule spanned by the arrows 
of the quiver, cf. \ref{MainProposition}.

\sus{The quiver relations}
First let us relabel the matrices $f_{y, w}$ using the 
enumeration of $y \in {\mathfrak S_3}$ as above.
Thus the matrix $f_{s_1s_2s_1, s_2s_1}$ will now be denoted
as $f_{(13)}$. Introduce the matrix $\widetilde d$ 
with the coefficients in $T(A_1)$ given by
\begin{equation}
\widetilde d =
\left(
\begin{matrix}
0 & (46) f_{(46)} & (56) f_{(56)} & 0 & 0 & 0 \\
(64) f_{(64)} & 0 & 0 & (24) f_{(24)} & (34) f_{(34)} & 0 \\
(65) f_{(65)} & 0 & 0 & (25) f_{(25)} & (35) f_{(35)} & 0 \\
0 & (42) f_{(42)} &  (52) f_{(52)} & 0 & 0 &  (12) f_{(12)} \\
0 & (43) f_{(43)} &  (53) f_{(53)} & 0 & 0 &  (13) f_{(13)} \\
0 & 0 & 0 & (21) f_{(21)} &  (31) f_{(31)} & 0
\end{matrix}
\right ).
\end{equation}
For example, 
$$
(42) f_{(42)} =
=\left(
\begin{matrix}
0 & 0 \\
(42) & 0 \\
-(42) & 0 \\
0 & (42)
\end{matrix}
\right).
$$
Some more notation: denote $f_{(ijk)} = f_{jk}\circ f_{ij}$, where
$\circ$ stands for matrix multiplication.

\sss{} Now, consider the $19\times 19$ matrix 
$\widetilde d^2$ as a matrix over $T(A_1)$.
$$  
{\widetilde d}^2 =
\left(
\begin{matrix}
g_{66} & 0 & 0 & g_{26} & g_{36} & 0 \\
0 & g_{44} & 
g_{54} & 0 & 0 & g_{14} \\
0 & g_{45} & g_{55}  & 0 & 0 & g_{15} \\
g_{62} & 0 & 0 & g_{22} & g_{32} & 0 \\
g_{63} & 0 & 0 & g_{23} & g_{33} & 0 \\
0 & g_{41} & g_{51}&  0 & 0 & g_{11}
\end{matrix}
\right).
$$
For example,
\begin{equation}\label{g23}
\begin{aligned}
g_{23} & = (243)f_{(243)} + (253)f_{(253)} + (213)f_{(213)} \\
& = 
\left(
\begin{matrix}
0 & 0 & 0 & 0 \\
(253) + (213) & 0 & 0 & 0 \\
(243) - (253) & 0 & 0 & 0 \\
0 & (253) + (213) & (243) + (213) & 0
\end{matrix}
\right).
\end{aligned}
\end{equation}
In general,
$$
\begin{aligned}
g_{66} & = (646)f_{(646)} + (656)f_{(656)} \\
g_{26} & = (246)f_{(246)} + (256)f_{(256)} \\
g_{36} & = (346)f_{(346)} + (356)f_{(356)} \\
g_{44} & = (464)f_{(464)} +  (424)f_{(424)} + (434)f_{(434)} \\
g_{54} & = (564)f_{(564)} + (524)f_{(524)} + (534)f_{(534)} \\
g_{45} & = (465)f_{(465)} + (425)f_{(425)} + (435)f_{(435)} \\
g_{55} & = (565)f_{(565)} + (525)f_{(525)} + (535)f_{(535)} \\
g_{14} & = (124) f_{(124)} + (134)f_{(134)} \\
g_{15} & = (125) f_{(125)} + (135)f_{(135)} \\
g_{62} & = (642)f_{(642)} + (652)f_{(652)} \\
g_{63} & = (643)f_{(643)} + (653)f_{(653)} \\
g_{22} & = (242)f_{(242)} + (252)f_{(252)} + (212)f_{(212)} \\
g_{32} & = (342)f_{(342)} + (352)f_{(352)} + (312)f_{(312)} \\ 
g_{23} & = (243)f_{(243)} + (253)f_{(253)} + (213)f_{(213)} \\
g_{33} & = (343)f_{(343)} + (353)f_{(353)} + (313)f_{(313)} \\
g_{41} & = (421) f_{(421)} + (431)f_{(431)} \\
g_{51} & = (521) f_{(521)} + (531)f_{(531)} \\
g_{11} & = (121) f_{(121)} + (131)f_{(131)}.
\end{aligned}
$$
Each of the matrices $g_{ij}$ above can be presented in the explicit 
form (\ref{g23}).

\sss{A basis of relators} Thus $\widetilde d^2$ is 
a matrix with elements in $A_1\otimes_k A_1$.
By section \ref{SectionQuiverd2}
the elements of $\widetilde d^2$ span
the submodule of relators $R \subset A_1\otimes_k A_1$.
One can choose a linearly independent subsystem (basis of $R$) 
from elements 
of $\widetilde d^2$ for example as follows:
\begin{equation}\label{relA2}
\begin{aligned}
\{ &  (121), (131), \\
& (242), (252) + (212), \\
& (353), (343) + (313), \\
& (243) + (213),(253) + (213), \\
& (352) + (312), (342) + (312), \\ 
& (124) + (134), (125) + (135), \\
&(246) + (256), (346) + (356), \\
& (421) + (431), (521) + (531),\\
& (464) - (424), (565) - (535), \\ 
& (425) + (435) + (465),
(524) + (534) + (564), \\
& (642) + (652),
(643) + (653) \}.
\end{aligned}
\end{equation}
There are $22$ relators as $\dim R = 22$. 

\sus{Open question}\label{openQuestionPMone} 
Notice that all the coefficients of basic paths of lengths two 
in (\ref{relA2})
are
$\pm 1$ as opposed to \cite{St03} where the coefficients    
are more general integers. 
It would be interesting to know whether
coefficients of basic paths of lengths two 
for the relators can be chosen to be $\pm 1$
in arbitrary type. 

\end{document}